\documentclass[reqno,12pt]{amsart}
\usepackage{epsf}
\newtheorem{theorem}{Theorem}

\newtheorem*{remark}{Remark}
\newtheorem*{cor}{Corollary}

\DeclareMathOperator{\Pf}{Pf}

\DeclareMathOperator{\stat}{stat}

\def\al{\alpha}

\def\ga{\gamma}
\def\la{\lambda}

\def\dsum{\displaystyle\sum}
\begin{document}
\thispagestyle{empty}

\newbox\Adr
\setbox\Adr\vbox{
\vskip.5cm
\centerline{Fakult\"at f\"ur Mathematik, Universit\"at Wien,}
\centerline{ Nordbergstra\ss{}e 15, A-1090 Wien, Austria.}
\centerline{E-mail: \footnotesize{\tt Theresia.Eisenkoelbl@univie.ac.at}}
}

\title[A Schur function identity for the $(-1)$--enumeration of SCPPs]{
\boldmath A Schur function identity related to the $(-1)$--enumeration of self-complementary plane partitions} 
\author[Theresia Eisenk\"olbl]{Theresia Eisenk\"olbl$^\ast$
\box\Adr
}
\thanks{ $^\ast$This research was supported
      by the Austrian Science Foundation FWF, grant S9607, in the
      framework of the National Research Network ``Analytic
      Combinatorics and Probabilistic Number Theory" and grant P19650
      ``Evaluation of Determinants and Pfaffians in Enumerative Combinatorics"}

\begin{abstract}
We give another proof for the $(-1)$--enumeration of
self-complementary plane partitions with at least one odd side-length 
by specializing a certain Schur
function identity. The proof is analogous to Stanley's proof for the
ordinary enumeration.  
In addition, we obtain enumerations of $180^\circ$-symmetric rhombus tilings of hexagons with a barrier of arbitrary length along the central line.
\end{abstract}

\maketitle

\begin{section}{Introduction} \label{introsec}

Plane partitions were first introduced by MacMahon (see
Figure~\ref{sceoofi} for an example and Section~\ref{defsec} for a
definition). 
He counted plane partitions contained in a given box \cite[Art.~429, 
proof in Art.~494]{MM} (see Eq.~\eqref{box}) and also investigated the
number of 
plane partitions with certain symmetries.

In \cite{MRR}, Mills, Robbins and Rumsey introduced additional
complementation symmetries giving six new combinations of symmetries
which led to more conjectures all of which were settled in the 1980's
and 90's (see \cite{Stan2,Ku2,An3,Stem4}). All these numbers can be
expressed as nice product formulas typically involving rising
factorials.

Many of these theorems come with $q$--analogs.
Recall that in a $q$--analog of an enumeration result for a set $M$,
each object is counted by a power of $q$, that is, one considers
$f_M(q)=\sum_{x\in M} q^{\stat(x)}$, where stat assigns an integer to
each object in $M$.

In the case of plane partitions, a typical example for
$\stat(x)$ is the number of little cubes in the plane partition $x$.
The closed forms for $f_{M}(q)$ now contain $q$--rising factorials instead
of rising factorials (see \cite{An1,An2,MRR2}). 

Interestingly, upon setting $q=-1$ in
the various $q$--analogs, one consistently obtains 
enumerations of other objects, usually with additional symmetry
constraints. This observation, dubbed the ``$(-1)$--phenomenon" has been
explained for many but not all cases 
by Stembridge (see \cite{Stem94a} and \cite{Stem94b}).

For the plane partitions with complementation symmetries, the
aforementioned $q$-weights either give trivial results or are not
well-defined. The one exception is the enumeration of self-complementary
plane partitions which was settled by Stanley \cite{Stan2} using a
Schur function identity. It gives rise to a $q$-analog via the
principal specialization (that is, setting $x_i=q^i$ in the Schur
functions, see Eq.~\eqref{factor}).

On the other hand, for all cases with complementation symmetries
Kuperberg defined a $\pm 1$-- weight \cite[p.26]{Ku} (that is, one considers
$\sum_{x\in M} w(x)$, where $w(x)=1$ or $w(x)=-1$). It has been
proved in Kuperberg's own paper and in \cite{min} and \cite{minsc}
that this sum admits a nice closed product formula in almost all cases. 

In accordance with the
``$(-1)$-phenomenon" mentioned above, the results coincide with the
enumerations of other similar objects, but in most cases an
explanation in the spirit of Stembridge's papers is missing.

The purpose of the present paper is the explanation of another curious
observation. One could expect to obtain a proof of the
$(-1)$--enumeration for self-complementary plane partitions from
\cite{minsc}
by setting 
$x_i=(-1)^i$ 
in the Schur function identity
Stanley uses in \cite{Stan2}. However, the $\pm 1$--weight on
individual plane partitions arising in this way is different. 
Yet, for self-complementary plane partitions with at least one odd
sidelength, setting $x_i=(-1)^i$ in Stanley's formulas yields exactly
the same formulas as in \cite{minsc}.

We will see that this mystery can be explained by a similar Schur
function identity which is actually a generalization 
of the one Stanley uses.
The three sides of the box containing the plane partitions originally
play a symmetrical role, but in the Schur function approach the
symmetry is broken arbitrarily. In \cite{Stan2} (see also the erratum),
this is done in a way to minimize complications. We will see
below that a less straightforward approach
produces the desired $\pm 1$--weight. As a by-product, we obtain
an additional
result about certain subclasses of self-complementary plane partitions
with a fixed line in the middle (see Figure~\ref{cut} and
Theorem~\ref{genenum}).
The Schur function identities given in Theorem~2 have already been obtained in
\cite[Cor 7.3]{IsOkTaZe} by Ishikawa, Okada, Tagawa and Zeng as a corollary to a very general Pfaffian
identity whose entries are products of determinants. In this paper, we
will give a different direct proof using the Littlewood-Richardson rule.

In Section~\ref{defsec}, we will review the necessary definitions and
properties of plane partitions and Schur functions. In
Section~\ref{schursec}, we state and 
prove the Schur function identity (see  Theorem~\ref{schur}).
Finally in
Section~\ref{ppsec}, we explain the connection with plane partitions
and draw the conclusions for the enumeration of self-complementary plane
partitions with a fixed line in the middle in Theorem~\ref{genenum} and
the $(-1)$--enumeration of
self-complementary plane partitions with at least one odd sidelength in
Theorem~\ref{minenum}.
The enumeration expressed as a product formula in
Theorem~\ref{genenum} can also be expressed as a Pfaffian using the
methods from \cite{minsc}, which immediately leads to the evaluation of
the Pfaffian stated in the subsequent 
corollary. In fact, the structure of the Pfaffian evaluation was
first observed by computer experiments and led backwards to the right
Schur function identity.

We remark that our results do not explain why the enumeration of self-complementary plane
partitions factors into terms corresponding to the enumeration of
ordinary ones, but an explanation should probably include the case of
self-complementary plane partitions with a fixed line in the middle
because the result has exactly the same structure.

\end{section}
\begin{section}{Definitions and basic properties} \label{defsec}

A plane partition $P$ can be defined as a finite set of points $(i,j,k)$
with integers $i,j,k > 0$ and if $(i,j,k) \in P$ and $1\le i'\le i$, 
$1\le j'\le j$, $1\le k'\le k$ then $(i',j',k')\in P$. We
interpret these points as midpoints of cubes and represent a plane
partition by stacks of cubes (see Figure~\ref{sceoofi}). If we have
$i\le a$, $j\le c$ and $k\le b$ for all cubes of the plane partition,
we say that the plane partition is contained in a box with sidelengths 
$a,b,c$.

\begin{figure}[t]
\begin{center}
\leavevmode
\epsfbox{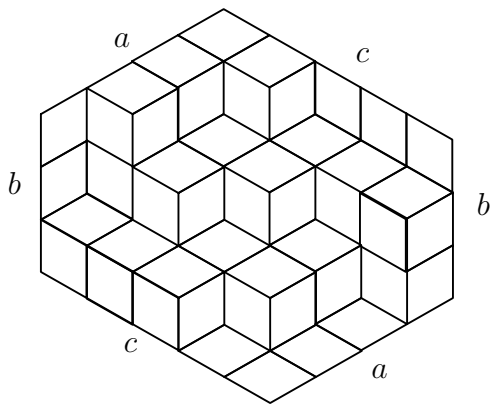}
\end{center}
\caption{A self--complementary plane partition}
\label{sceoofi}
\end{figure} 

Another way to represent a plane partition is by writing down the
number of cubes above a certain place in the $xy$--plane giving an array of
numbers with weakly decreasing rows and columns. The plane
partition of Figure~1 corresponds to the array 
\begin{equation} \label{array}
\begin{matrix}
3& 3& 2& 2& 2\\
3& 2& 2& 1& 0\\
3& 2& 1& 1& 0\\
1& 1& 1& 0& 0
\end{matrix}
\end{equation}

This representation gives the connection to semi-standard
tableaux and Schur functions
which was used in Stanley's proof of the ordinary enumeration of
self--complemen\-tary plane partitions (see Section~\ref{ppsec}).

A plane partition $P$ contained in an $a\times b\times c$--box
is called
{\em self--complementary}
if  whenever $(i,j,k) \in P$ then $(a+1-i,c+1-j,b+1-k) \notin P$ for 
$1\le i\le a$, $1\le j\le c$, $1\le k\le b$ (see Figure~\ref{sceoofi}).

The corresponding array of numbers has the property that an entry and
the corresponding entry in the array rotated by $180^\circ$ add up to
$b$ (see the array in \eqref{array} with $b=3$).

\begin{figure}
\begin{center}
\leavevmode
\epsfbox{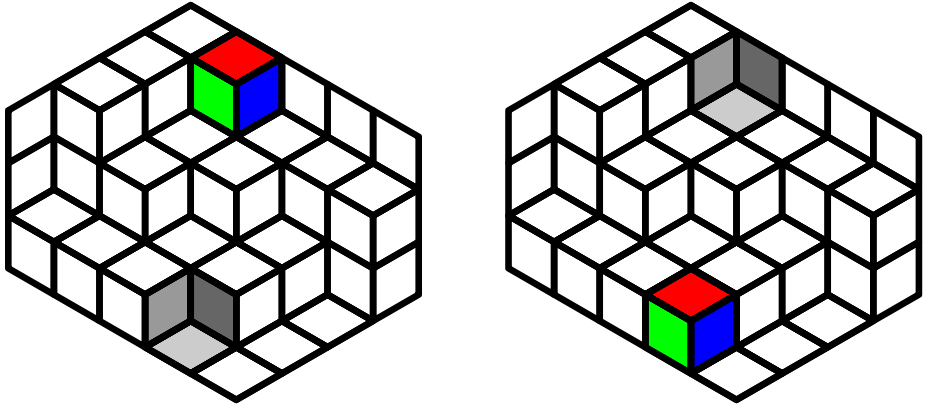}
\end{center}
\caption{Two plane partitions with different $\pm 1$--weight}
\label{oppfi}
\end{figure}

Now, we define the $\pm 1$-weight for self-complementary plane partitions.
The weight
changes sign if we remove one cube and add the opposite one
(see Figure~\ref{oppfi}).
The half-full plane partition has weight 1 (see
Figure~\ref{sceoonormfi}). Since every self-complementary plane
partition can be reached from the half-full plane partition by
moving cubes as described above, this defines the weight uniquely.

\begin{figure}
\begin{center}
\leavevmode
\epsfbox{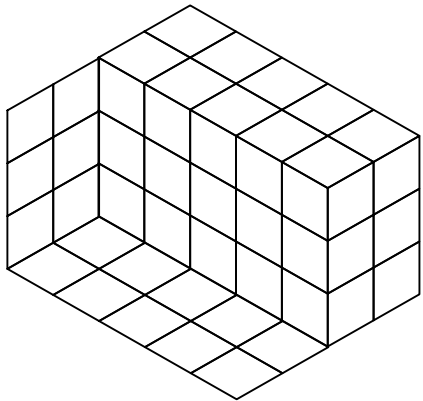}
\end{center}
\caption{A plane partition of weight 1.}
\label{sceoonormfi}
\end{figure}

A partition $\la$ is a sequence of integers $\la_1 \ge \dots \ge
\la_l$ which can be represented by a Ferrers diagram as an array of
squares with $\la_i$ squares in row $i$. For example, the partition
$(4,2,1)$ is represented by the following Ferrers diagram.

{
\begin{center}
\epsfysize=1cm
\epsfbox{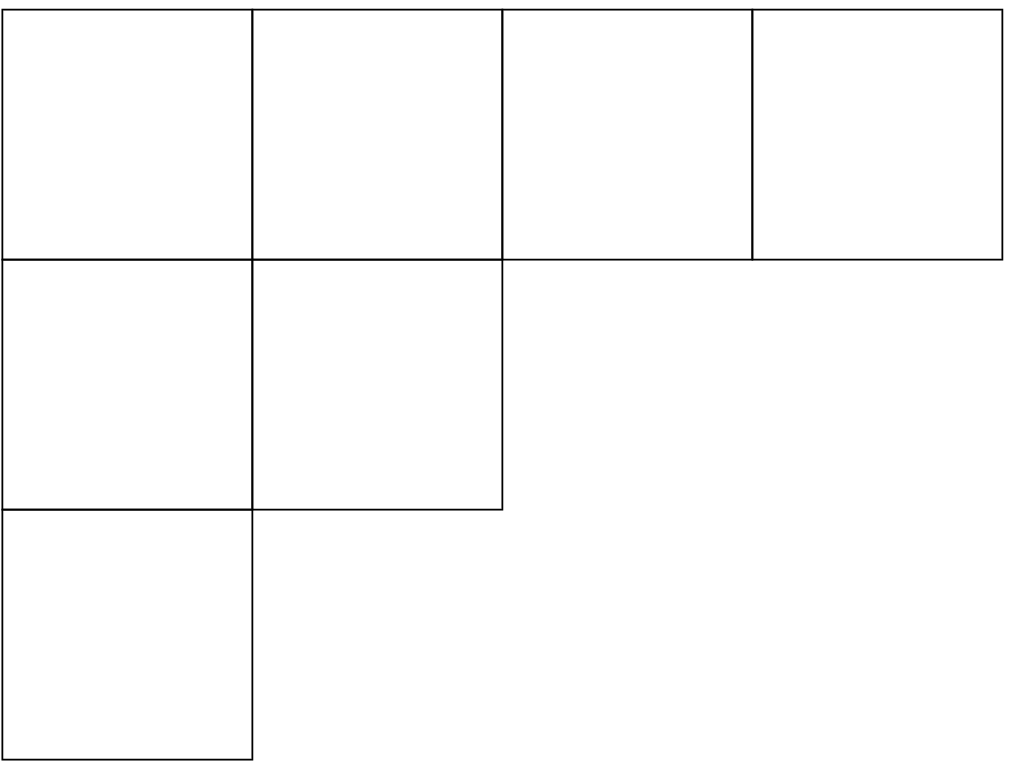}
\end{center}
}

A partition $\mu$ is contained in a partition $\la$ if $\mu_i\le
\la_i$ for all $i$. $\la/\mu$ denotes the set of squares in the Ferrers
diagram of $\la$ that are not in the Ferrers diagram of $\mu$. If this
set of squares does not contain two squares in the same column, it is
called a horizontal strip. By $|\la/\mu|$, we denote the number of
squares in $\la/\mu$.

A semistandard tableau of shape $\la$ (or $\la/\mu$) is a filling of
the squares of 
the Ferrers diagram of $\la$ (or $\la/\mu$) with positive integers which is
increasing along the rows and strictly increasing along the
columns. If there are $\nu_1$ entries $1$, $\nu_2$ entries $2$, etc.,
then $\nu=(\nu_1,\dots,\nu_k)$ is called the content of the
tableau $T$.

A Schur function is the generating function of semistandard tableaux of
a given shape, namely $s_{\la}(x_1,\dots, x_n) = \sum_{T} x^T$ where $T$ runs over
the semistandard tableaux of shape $\la$ and 
$x^T=x_1^{\text{\# 1's in $T$}}x_2^{\text{\# 2's in $T$}}\dots$.

We will use the following results about plane partitions. The first
result was obtained by MacMahon in
\cite[Art.~429,$x\to1$, proof in Art.~494]{MM}:

\smallskip
{\em The number of all plane partitions contained in a box 
with sidelengths $a,b,c$ equals
\begin{equation} \label{box}
B(a,b,c)= \prod _{i=1} ^{a}\prod _{j=1} ^{b}\prod _{k=1} ^{c}\frac {i+j+k-1}
{i+j+k-2}= \prod _{i=1} ^{a}{\frac {(c+i)_b} {(i)_b}}, 
\end{equation}
where $(a)_n:=a(a+1)(a+2)\dots (a+n-1)$ is the rising factorial.}

Stanley's proof of the ordinary enumeration of self--complementary
plane partitions gives the following result:

\begin{theorem}[Stanley \cite{Stan2}] \label{th:stanley}
The number $SC(a,b,c)$ of self--complementary plane partitions
contained in a box with sidelengths $a,b,c$ can be expressed in terms
of $B(a,b,c)$ in the following way:

\begin{align*}
B\left(\tfrac a2,\tfrac b2, \tfrac c2\right)^2 \quad &\text{for $a,b,c$ even,}\\
B\left(\tfrac a2,\tfrac{b+1}2,\tfrac{c-1}2\right) B\left(\tfrac
a2,\tfrac{b-1}2,\tfrac{c+1}2\right) \quad &\text{for $a$ even and $b$, $c$
  odd,}\\
B\left(\tfrac {a+1}2,\tfrac{b}2,\tfrac{c}2\right) B\left(\tfrac
{a-1}2,\tfrac{b}2,\tfrac{c}2\right)
\quad &\text{for $a$ odd and $b$, $c$
  even.}
\end{align*}
\end{theorem}

The result is obtained as the case $x_i=1$ of the expansions in Schur
functions of the products 
\begin{align*}
&s_{s^r}(x_1,x_2,\dots,x_{t+r})^2 \quad &&\text{for the $2r\times 2s
  \times 2t$--box,}\\ 
&s_{s^r}(x_1,x_2,\dots,x_{t+r})s_{(s+1)^r}(x_1,x_2,\dots,x_{t+r}) \quad
&&\text{for the $2r \times (2s+1) \times 2t$--box,}\\  
&s_{s^{r+1}}(x_1,x_2,\dots,x_{t+r+1})s_{s^r}(x_1,x_2,\dots,x_{t+r+1}) \,\,
&&\text{for the $(2r+1) \times 2s \times (2t+1)$--box.}
\end{align*}

Note that for $x_i=q^i$, $i=1,\dots, n$, 
a Schur function corresponding to a rectangular partition
can be expressed by a product formula as a special case of Stanley's
hook-content formula (\cite{Stan3}):  
\begin{equation} \label{factor}
s_{\ga^\al}(q,q^2,\dots,q^n)=q^{\ga \al(\al+1)/2} s_{\ga^\al}(1, q,
\dots,q^{n-1})= 
q^{\ga \al (\al+1)/2} \prod_{i=1}^{\al} \prod_{k=0}^{\ga-1} \frac
{1-q^{i+n-\al+k}} {1-q^{i+k}}.
\end{equation}

Furthermore, it is straightforward to obtain from this equation that
$s_{c^a}(1,1,\dots,1)$ (with $a+b$ arguments of 1) is $B(a,b,c)$
and $s_{c^a}(1,-1,1,\dots,¸(-1)^{a+b-1})$ (with $a+b$ arguments in the
Schur function) is
$SC(a,b,c)$ (see Section~2 of \cite{Stem94a}).

\end{section}

\begin{section}{The Schur function identity} \label{schursec}

In this section, we state and prove a Schur function identity which
implies a closed form for the $(-1)$-enumeration of self-complementary
plane partitions (see Theorem~\ref{minenum}).

\begin{theorem}[Ishikawa, Okada, Tagawa, Zeng \cite{IsOkTaZe}] \label{schur}

Let $\ga_1 \ge \ga_2$. Then we have the following two identities:

\begin{equation} \label{schurid1}
\begin{aligned} 
s_{(\ga_1^\al)}(x_1,\dots, x_{n})s_{(\ga_2^\al)}(x_1,\dots,x_{n+1}) 
\hskip-5cm&\\
&= \sum _{\la \subseteq (\ga_2^\al)}\sum _{\substack{\pi \subseteq \la
    \\ \la
  / \pi \text{ horiz. strip 
}}} x_{n+1}^{|\la/\pi|}\cdot
s_{(\ga_1+\ga_2-\la_{\al},\ga_1+\ga_2-\la_{\al-1},\dots,
  \ga_1+\ga_2-\la_{1},
\pi_1,\dots,\pi_{\al} )}(x_1,\dots, x_n)\\
&= \sum _{\la \subseteq (\ga_2^\al)}\sum _{T} x^T,
\end{aligned}
\end{equation}
where $T$ is a semistandard tableau of the following shape:

{ 
\epsfysize=5cm
\epsfbox{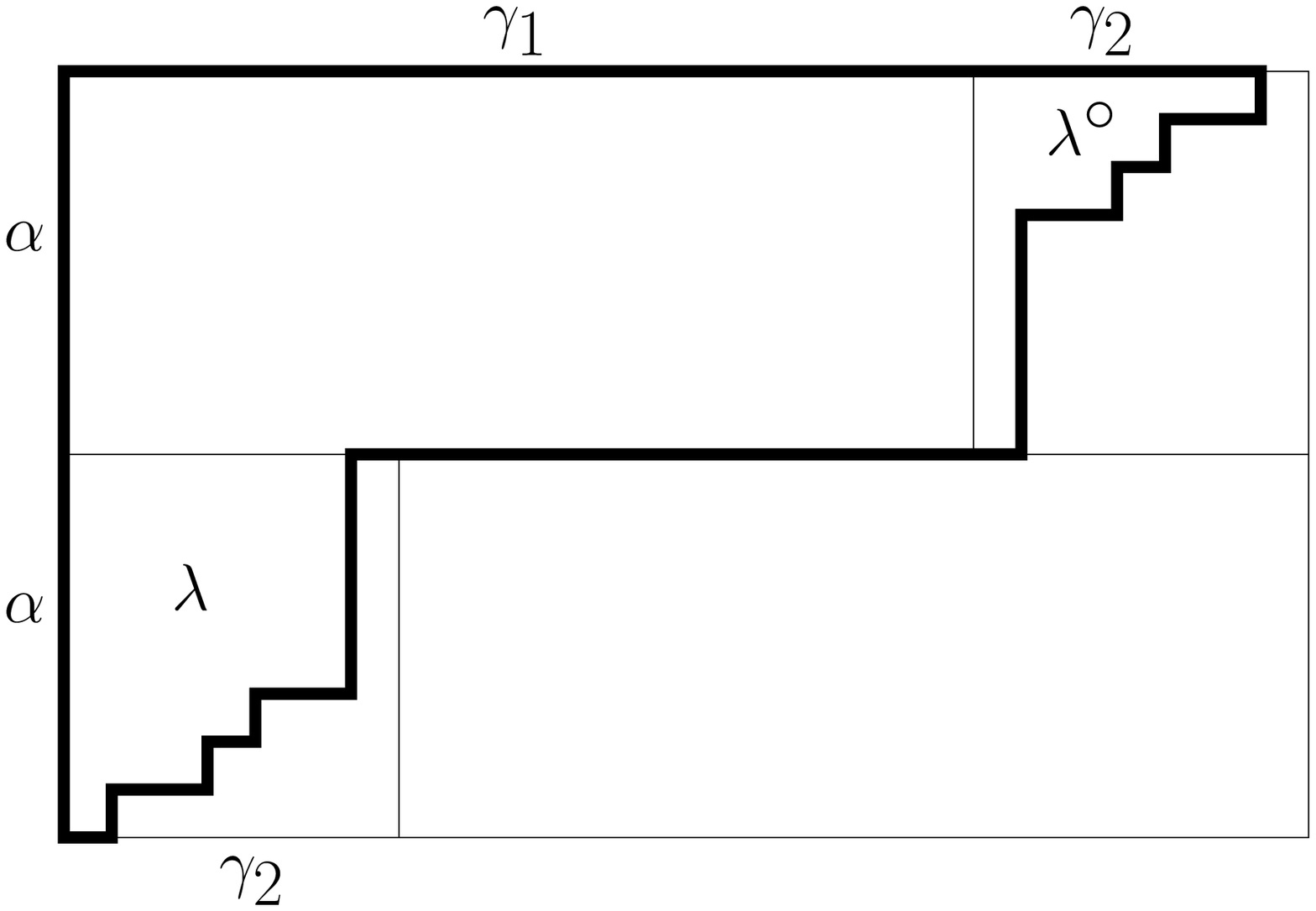}
}

Here, entries $n+1$ can only occur in the region of shape
$\la$ and all other entries are smaller. $\la^\circ$ is the shape
obtained by taking the complement of $\la$ in the rectangle
$\ga_2^{\al}$ 
and rotating by $180^\circ$. $x^T$ is short-hand notation for
the monomial $x_1^{\text{\# 1's in $T$}}x_2^{\text{\# 2's in $T$}}\dots$.


\begin{equation} \label{schurid2}
\begin{aligned} 
s_{(\ga_1^\al)}(x_1,\dots, x_n)s_{(\ga_2^{\al+1})}(x_1,\dots,
x_{n+1})\hskip-6cm&\\ 
&=\sum _{\substack{\la \subseteq (\ga_2^{\al+1})\\ \la_1=\ga_2}}
\sum _{\substack{\pi \subseteq \la\\ \la / \pi \text{
      horiz. strip
}}}  x_{n+1}^{|\la/\pi|} \cdot
s_{(\ga_1+\ga_2-\la_{\al+1},\ga_1+\ga_2-\la_{\al},\dots,
  \ga_1+\ga_2-\la_{2}, 
\pi_1,\dots,\pi_{\al+1} )}(x_1,\dots,x_n)\\
&= \sum _{\substack{\la \subseteq (\ga_2^{\al+1})\\ \la_1=\ga_2}}\sum _{T} x^T,
\end{aligned}
\end{equation}
where $T$ is a semistandard tableau of the following shape:

{ 
\epsfysize=5cm
\epsfbox{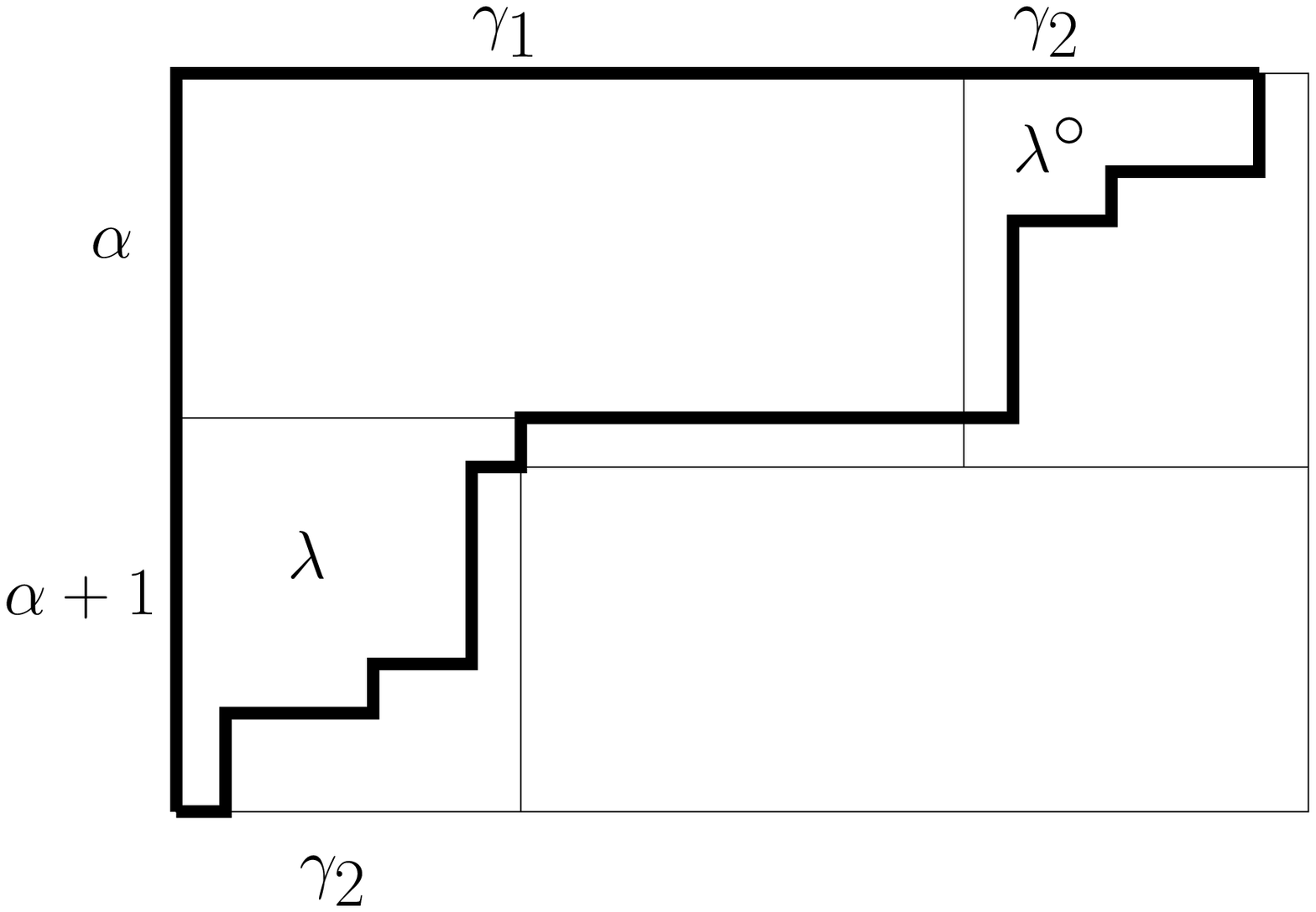}
}

Entries $n+1$ can again only occur in the region of shape $\la$.The
first row of $\la$ is $\ga_2$, and $\la^\circ$ is the rotated
complement of $\la$ in the rectangle $\ga_2^{\al+1}$.

\end{theorem}

\begin{remark}
The case $x_{n+1}=0$ and $\ga_1=\ga_2$ of \eqref{schurid1} and
\eqref{schurid2} gives the identities used by Stanley. Note also that
both sides of the identities are zero if $n < \al$.
\end{remark}

\begin{proof}


We want to expand the product on the left-hand side of the first
equation.
 The tableaux in
the second factor can contain entries $n+1$, but only in the last
row. Therefore, we can split the term into summands in the following
way:

\begin{multline*}
s_{(\ga_1^\al)}(x_1,\dots, x_n)\cdot s_{(\ga_2^{\al})}(x_1,\dots,
x_{n+1}) \\
=\sum _{k=0}^{\ga_2} x_{n+1}^k \cdot s_{(\ga_1^\al)}(x_1,\dots,x_n) \cdot
s_{(\ga_2^{\al-1},\ga_2-k)}(x_1,\dots, x_n).
\end{multline*}

This product of two Schur functions can be expanded by the
Littlewood--Richardson rule (see \cite{LR}), i.e. $s_\mu(x_1,\dots,x_n)\cdot
s_\nu(x_1,\dots, x_n)
=\sum_{\rho}
c_{\mu\nu}^{\rho} s_\rho (x_1,\dots, x_n)$, where $c_{\mu\nu}^{\rho}$ is the number of
semistandard tableaux of shape $\rho / \mu$ and content $\nu$ such
that reading each 
row from right to left starting with the first row we always have $\#
1's \ge \# 2's \ge \# 3's \ge \dots$.

We apply the Littlewood-Richardson rule with $\mu=(\ga_1^\al)$ and
$\nu=(\ga_2^{\al-1},\ga_2-k)$. 
There are not many possibilities for the shape $\rho$ which has to
contain the shape $\mu$ and must have $|\mu|+|\nu|$ boxes.

{ 
\epsfysize=5cm
\epsfbox{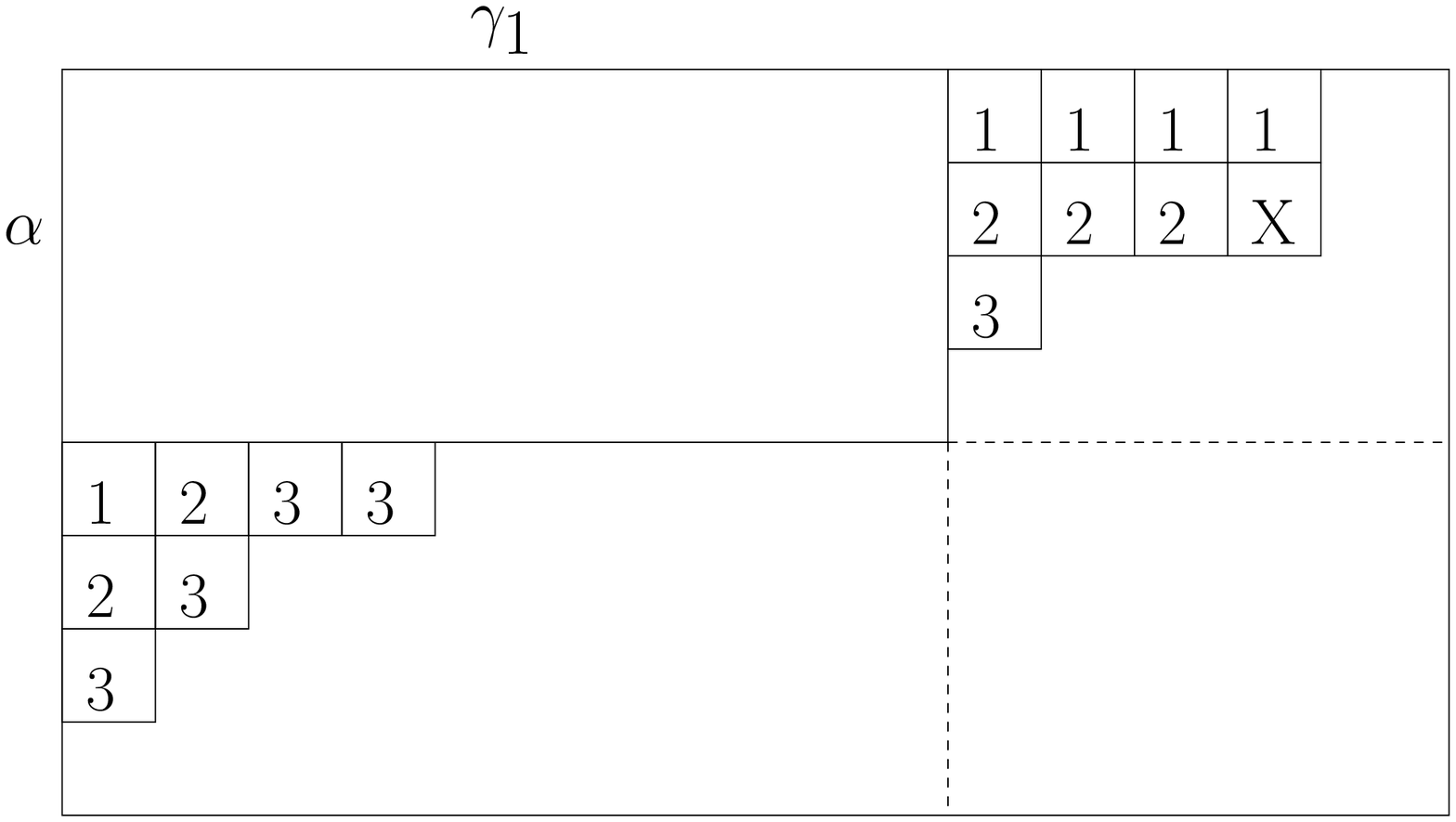}
}

The left upper corner of this picture is the shape $\mu=(\ga_1^\al)$
which has been taken out of $\rho$. The 
right lower corner is also empty since the partition 
$\nu=(\ga_2^{\al-1},\ga_2-k)$
encoding the 
content tells us that the highest possible entry in the
Littlewood-Richardson tableau of shape $\rho/\mu$ is $\al$.

It turns out that row $i$ in the right upper block consists only of
entries $i$. Consider for example the position $X$ in the
picture. It cannot contain an entry $3$ (or higher) because the word
$11113222\dots$ would have a $3$ before the first $2$.
Since the number of 1's is given by $\nu_1=\ga_2$, we can see that the
right upper corner is inside the rectangle $\al \times \ga_2$, and
there is exactly one admissible filling of each partition $\la^\circ$
inside this rectangle. 
 
Now, how many possibilities are there to place the remaining
$\la_{\al}$ entries 1, $\la_{\al-1}$ entries 2, $\dots$, $\la_2$
entries $\al-1$ and
$\la_1-k$ entries $\al$ into the left lower corner?

We claim that for the entries from 1 to $\al-1$ there is only the
following possibility: Write all 1's in the first row, then write a 2
under each 1 and place the remaining 2's in the first row, then write
a 3 under each 2 and place the remaining 3's in the first row, and so
on up to $\al-1$.

Suppose we have already filled in the entries from $1$ to $i$ in this
way (in particular, every column ends with $i$ at this point and we
can write entries $i+1$ either in the first row or underneath an entry
$i$).
Now assume that we would not write $i+1$ under each of the $\la_{\al-i+1}$
entries $i$.

 Then we read the Littlewood-Richardson word up to and including
the entries $i+1$ 
of the first row of the left lower corner (if there is no entry $i+1$
in the first row, we just read the right upper corner).
In this word, we have
$\la^\circ_i=\ga_2-\la_{\al-i+1}$ 
entries $i$ and at least $\ga_2-\la_{\al-i+1}+1$ entries $i+1$ (all but those
in the lower left corner under an entry $i$) which clearly contradicts
the Littlewood-Richardson condition.

{ 
\epsfysize=5cm
\epsfbox{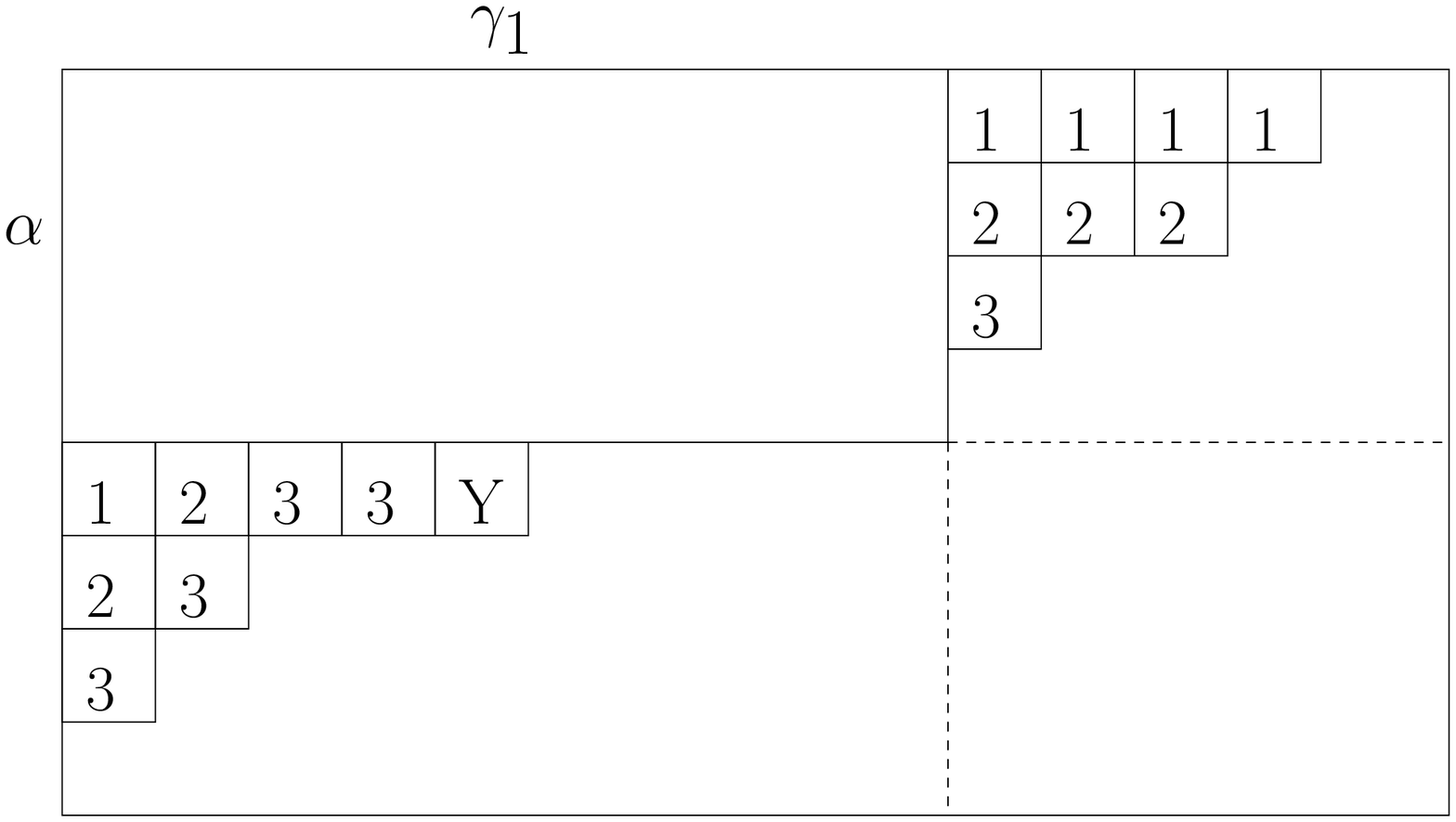}
}

Consider for
example the position $Y$ in the picture. It cannot contain an entry
$3$ because we would read four 3's before the fourth 2. Therefore, we
have to write a 3 underneath each 2.

Note that this argument also proves that the first row of the
partition in the left lower corner cannot be bigger than $\ga_2$.

In the case $k=0$, the same argument holds for the entries $\al$ and
the shape in the left lower corner is $\la$ with a unique filling.
If $k\not=0$, then we certainly cannot have more $\al$'s in the first
row, so we simply have to remove some of the $\al$'s which means that
we remove a horizontal strip from the shape $\la$. Since the filling
is unique for each choice of the horizontal strip all
the coefficients $c_{\mu\nu}^{\rho}$ are 1.

This leads immediately to the first equality in identity \eqref{schurid1}.
We get the second equality by filling the missing
horizontal strip with entries $n+1$ in the resulting tableaux of shape
$\rho$ occuring in the Schur function. This gives tableaux of the shape
$T$ in the figure where entries $n+1$ occur only in a certain region.

The second identity is proved analogously. We just have to make sure
that there can be no $\al+1$ in $\la^\circ$ which follows from the
fact that the first row of the shape in the left lower corner cannot
be bigger than $\ga_2$ and $\ga_2\le \ga_1$.
\end{proof}
\end{section}

\begin{section}{Connection to self-complementary plane partitions} \label{ppsec}

In this section we will explain how to get the $(-1)$--enumeration for
self-complemen\-tary plane partitions from Theorem~\ref{schur} for the
case of a box with at least one odd side and why this does not work
with the Schur function identity from Stanley's proof.

Furthermore, we will explain what kind of self-complementary plane
partitions are counted by the general case $\ga_1\not=\ga_2$.

We will prove the following two theorems.

\begin{theorem} \label{genenum}
The enumeration of self-complementary plane partitions in an $a\times
b\times (c_1+c_2)/2$--box where the corresponding rhombus tiling
contains a fixed ``middle line" of length $(c_1-c_2)/2$ parallel to
$(c_1+c_2)/2$  (see Figure~\ref{cut}) equals 

\begin{align*}
B\left(\tfrac a2,\tfrac b2, \tfrac {c_1}2\right) B\left(\tfrac a2,\tfrac b2, \tfrac
{c_2}2\right)
 \quad &\text{for $a,b,c_1,c_2$ even,}\\
B\left(\tfrac {a+1}2,\tfrac{b}2,\tfrac{c_1}2\right) B\left(\tfrac
{a-1}2,\tfrac{b}2,\tfrac{c_2}2\right)
\quad &\text{for $a$ odd and $b$, $c_1$,$c_2$
  even,}\\
 B\left(\tfrac {a-1}2,\tfrac{b+1}2,\tfrac{c_1}2\right) B\left(\tfrac
{a+1}2,\tfrac{b-1}2,\tfrac{c_2}2\right) \quad &\text{for $a,b$ odd and $c_1$, $c_2$
  even.}
\end{align*}
where $B(a,b,c)$ is given in \eqref{box}.

In all cases the ``middle line" lies at equal distance from the border
of the hexagon in the direction of the line and orthogonal to the line
(see Figure~\ref{cut}).

The ``middle line" is just a line of length $\frac {c_1-c_2} {2}$ if
$a,b$ is even, it is a strip of $\frac {c_1-c_2} {2}$ rhombi if
$a$ is even and $b$ is odd and it is a line of length $\frac {c_1
-c_2} {2}-1$ with two attached triangles if $a$ is odd
and $b$ is even.
\end{theorem}

In terms of the corresponding integer arrays $a\times \frac {c_1+c_2}
{2}$ with weakly decreasing rows and columns where opposite entries
add up to $b$, we have the condition that the entry in row
$a/2$ and column $c_1/2$ has to be at least $b/2$ in the first case.
In the second case, the condition is that the entries in row $(a+1)/2$ and
columns $c_2/2+1,\dots, c_1/2$ are filled
with entries $b/2$.
In the third case, the places  in row $(a+1)/2$ and
columns $c_2/2+1,\dots, c_1/2$ are empty, the entries above them are
at least $(b-1)/2$, the entries below them are at most $(b+1)/2$ and
the entry to the left of them is at least $(b+1)/2$.

\begin{theorem}\label{minenum}
The $(-1)$--enumeration of self-complementary plane partitions in an
$a\times b \times c$--box equals up to sign
\begin{align*}
&SC\left(\tfrac a2,\tfrac{b+1}2,\tfrac{c-1}2\right) SC\left(\tfrac
a2,\tfrac{b-1}2,\tfrac{c+1}2\right) \quad &&\text{for $a$ even and $b$, $c$
  odd}\\
&SC\left(\tfrac {a+1}2,\tfrac{b}2,\tfrac{c}2\right) SC\left(\tfrac
{a-1}2,\tfrac{b}2,\tfrac{c}2\right)
\quad &&\text{for $a$ odd and $b$, $c$
  even}
\end{align*}
where $SC(a,b,c)$ is given by Theorem~\ref{th:stanley}.
\end{theorem}

\begin{remark}
In the case of three even sidelengths, the $(-1)$--enumeration of
self-comple\-men\-tary plane partitions is $B(a/2,b/2,c/2)$ (see
\cite{min}), so we cannot expect to prove it 
with the help of an identity involving the product of two Schur functions.
\end{remark}

Similarly to the proof of the $(-1)$--enumeration in \cite{minsc}, 
the enumerations in Theorem~\ref{genenum}
can be expressed by Pfaffians via families of
non-intersecting lattice paths. Therefore, we have the following
corollary.

\begin{cor} \label{coro}
The following Pfaffians evaluate to the respective expressions in
Theorem~\ref{genenum}.

\begin{multline*} 
\Pf_{1\le i,j\le a} \left(\sum _{k=1}^{\frac {a+b} {2}}
\left( \binom{(b+c_1)/2}{b+i-k}\binom{(b+c_2)/2}{j+k-a-1}-
  \binom{(b+c_1)/2}{b+j-k}\binom{(b+c_2)/2}{i+k-a-1} \right) \right)
\\
\shoveright {\text{for $a,b,c_1,c_2$ even}}\\
\shoveleft (-1)^{\frac{a-1}2} \Pf_{1\le i,j \le a+1} 
\left( \begin{array}{c|c} 
\dsum _{k=1} ^{\frac {a+b-1} {2}} (\tbinom{\frac{b+c_1}2}{b+i-k}
\tbinom{\frac{b+c_2}2}{j+k-a-1}
 - \tbinom{\frac{b+c_1}2}{b+j-k} \tbinom{\frac{b+c_2}2}{i+k-a-1}) & 
\tbinom{\frac{b+c_1}2}{b+i-\frac {a+b+1} {2}}\\
\hline
-\tbinom{\frac{b+c_1}2}{b+j-\frac{a+b+1} 2 } & 0  
\end{array}\right)\\
\shoveright {\quad \text{for $a$ odd and
  $b,c_1,c_2$ even,}}\\
\shoveleft (-1)^{\frac{a-1}2} \Pf_{1\le i,j \le a+1} 
\left( \begin{array}{c|c} 
\dsum _{k=1} ^{\frac {a+b} {2}} (\tbinom{\frac{b+c_2-1}2}{b+i-k}
\tbinom{\frac{b+c_1+1}2}{j+k-a-1}
 - \tbinom{\frac{b+c_2-1}2}{b+j-k} \tbinom{\frac{b+c_1+1}2}{i+k-a-1}) & 
-\tbinom{\frac{b+c_2-1}2}{b+i-\frac {a+b} {2}-1}\\
\hline
\tbinom{\frac{b+c_2-1}2}{b+j-\frac{a+b} 2 -1} & 0  
\end{array}\right)\\
\quad \text{for $a,b$ odd and
  $c_1,c_2$ even,
}
\end{multline*}
The second and the third matrix should be read as $a\times a$--matrices with an
extra row and column.
\end{cor}


\begin{proof}[Proof of Theorem~\ref{genenum}]
We will see that the self-complementary plane partitions with this
additional constraint are in bijection with the semi-standard tableaux
appearing in the Schur function identities \eqref{schurid1} and
\eqref{schurid2}.

Let us start with the case $a,b,c_1,c_2$ even.

As stated in Section~\ref{defsec}, self-complementary plane partitions can
be represented by rectangular $a\times (c_1+c_2)/2$--arrays of positive integers with
decreasing rows and columns with the additional condition that
entries related by a $180^\circ$-rotation add up to $b$.
For example, the self-complementary plane partition in
Figure~\ref{cut}a
corresponds to the array

$$\begin{matrix}
4& 3& 3& 3& 3& 3& 2& 1\\
\cline{4-5}
3& 2& 1& 1& 1& 1& 1& 0
\end{matrix}
$$
The fixed bold line in the picture becomes the constraint that
array entries above the horizontal middle line of length $(c_1-c_2)/2$ are
greater or equal to $b/2$.

To switch from decreasing to increasing and to make the columns {\em
  strictly} increasing, we rotate the array by $180^\circ$ degree and
add $i$ to row $i$ and get

\begin{equation}\label{counter}\begin{matrix}
1& 2& 2& 2& 2& 2& 3& 4\\
\cline{4-5}
3& 4& 5& 5& 5& 5& 5& 6
\end{matrix}
\end{equation}

Now, we have a semi-standard tableau where entries related by a
$180^\circ$-rotation add up to $a+b+1$ 
which happens to be odd in this case. The entries above the middle
line 
of length $(c_1-c_2)/2$ are now smaller or equal to $(a+b)/2$.

If we now look only at the subtableau consisting of the entries $\le
(a+b+1)/2$, we get shapes $\rho$ with $\rho=\rho^\circ$ who are above
the middle line. This is exactly the shape of the tableau $T$ in
\eqref{schurid1}, but without the condition on the entries $n+1$.

Therefore, we can count these plane partitions by setting
$\ga_1=c_1/2$, $\ga_2=c_2/2$, $\al=a/2$, $n=(a+b)/2$, $x_1=x_2=\dots=x_n=1$ and
$x_{n+1}=0$.

(The case $\ga_1=\ga_2$  is an instance of the original proof of Stanley.)

The closed form for this enumeration follows at once from the
left-hand side of \eqref{schurid1} and 
identity~\eqref{factor} with $q=1$.

Similarly, we get the second case $a$ odd and $b,c_1,c_2$ even by setting
$x_1=x_2=\dots=x_{n+1}=1$, $n=(a+b-1)/2$, $\ga_1=c_1/2$,
$\ga_2=c_2/2$ and $\al=(a-1)/2$ in the  
second part of Theorem~\ref{schur} with an example shown in
Figure~\ref{cut}b.

The corresponding tableau is

$$\begin{matrix}
1& 1& 2& 2& 2& 3& 3& 3\\
2& 4& 4& {\fbox {4 }}& {\fbox{ 4}}& 4& 4& 6\\
5& 5& 5& 6& 6& 6& 7& 7\\
\end{matrix}
$$

After substracting $i$ from row $i$ and rotating, we obtain an array
with weakly decreasing rows and columns where the boxed entries become
$b/2$. This translates to a strip of $(c_1-c_2)/2$ rhombi in the
middle of the hexagon (see Figure~\ref{cut}b).

The last case $a,b$ odd and $c_1,c_2$ even is obtained by setting
$x_1=x_2=\dots=x_{n}=1$, $x_{n+1}=0$, $n=(a+b)/2$, $\ga_1=c_1/2$,
$\ga_2=c_2/2$ and $\al=(a-1)/2$ in the  
second part of Theorem~\ref{schur} with an example shown in
Figure~\ref{cut}c.

Since there are now no entries $n+1$, we can complete the
semi-standard tableau to a rectangle where corresponding entries add
up to $2n+1$. This means that there are no integers that properly fit
in the center part of the middle row and we can think of the entries
as $n+1/2$ to preserve all monotony restrictions.

In the example, the corresponding ``tableau'' with strictly increasing
columns is 

$$\begin{matrix}
1& 2& 3& 3& 4\\
2& {\fbox {3.5}}& {\fbox {3.5}}& {\fbox {3.5}}&  5\\
3& 4& 4& 5& 6
\end{matrix}
$$

After substracting $i$ from row $i$ and rotating, 
we obtain an array with weakly decreasing rows and columns where the
entries in row $(a-1)/2$ and columns $c_2/2+1, c_2/2+2, \dots, c_1/2$
are at least $(b-1)/2$ and the corresponding entries in row $(a+3)/2$
are at most $(b+1)/2$. Furthermore, the entry in row $(a+1)/2$ and
column $c_2/2$ is at least $(b+1)/2$.

If we represent this by stacks of cubes and treat the non-integer
entries as empty, we
obtain a bijection to rhombus tilings with $180^\circ$-rotational symmetry
of hexagons with sides $a,b,(c_1+c_2)/2$ and a middle line of length
$(c_1-c_2)/2-1$ with two triangles attached to the end (see Figure~\ref{cut}).

\end{proof}

\begin{remark}
Obviously, we could also look at the specialization involving
$x_{n+1}$ in the second part, which works, but gives plane partitions
which are a reflection of the ones obtained in the second case above.
The $c_1,c_2$ are always even, but this is no restriction on the
parity of $(c_1+c_2)/2$, so indeed all possible sidelengths of
hexagons are treated.
\end{remark}

\begin{proof}[Proof of Theorem~\ref{minenum}]

We want to set $x_{i}=(-1)^i$ in Theorem~\ref{schur} to obtain our result.

Let us look again at the example of a tableau corresponding to a
self-complementary plane partition.



$$
\begin{matrix}
1& 2& 2& 2& 2& 2& 3& 4\\
3& 4& 5& 5& 5& 5& 5& 6
\end{matrix}
$$

Clearly, the action on the tableau
corresponding to removing and adding a cube in a plane partition as
shown in Figure~\ref{oppfi} consists of adding one to an entry and
subtracting one of the entry related by a $180^\circ$-rotation.
If this action changes two entries 2 and 5  to  3 and 4, the
subtableau of entries $\le 3$ loses an entry 2 and gains an entry
3. Therefore,  the contribution of
these entries to the terms in the Schur function changes from $x_2$ to
$x_3$ which gives the desired sign change for $x_i=(-1)^i$.

On the other hand, if we exchange two entries 3 and 4, the weight
obviously remains unchanged, so in this case, we do not obtain the
desired $\pm 1$-weight.



In this example, the maximal possible
entry $a+b$ in the array is even and the only case leading to problems are
two entries $(a+b)/2$ and $(a+b)/2+1$ swapping place leaving the
weight in the Schur function unchanged while the $\pm 1$--weight changes.

Since the parameters $a,b,c$ played a symmetrical in the definition of
self-complemen\-tary and of the $\pm 1$--weight, we can assume that $a+b$
is odd if not all of $a,b,c$ are even. 

Now, the only problematic case involves two entries $(a+b+1)/2$
changing to $(a+b+1)/2+1$ and $(a+b+1)/2-1$. It is clear that exactly one
of the entries $(a+b+1)/2$ must be in the left lower region where it
contributes to the weight, so one factor in the weight changes from
 $x_{(a+b+1)/2}$ to $x_{(a+b+1)/2-1}$ which coincides with the change of
 the $\pm 1$--weight for $x_i=(-1)^i$.

So, we assume that $a$ and $b$ have different parity and it only
remains to set
$\al=a/2$, $\ga_1=(c+1)/2$, $\ga_2=(c-1)/2$ and $n=(a+b-1)/2$  in
\eqref{schurid1} for 
the case $a$ even and $b,c$ odd  and
$\al=(a-1)/2$, $\ga_1=\ga_2=c/2$ and $n=(a+b-1)/2$  in \eqref{schurid2} for
the case $a$ odd and $b,c$ even.
(Note that the condition imposed by the middle-line of length 1 is
automatically fulfilled in a self-complementary plane partition and
therefore does not change the enumeration.)

The left-hand side in both cases lead to the case $q=-1$ of
\eqref{factor} which gives the desired result.
\end{proof}

\begin{proof}[Proof of the Corollary]
The original proof of Theorem~\ref{minenum} used a bijection between
plane partitions and families of non-intersecting lattice paths which
gives a Pfaffian for both the weighted and the unweighted case which
can then be evaluated.  

Since this is the case $c_1=c_2$ in Theorem~\ref{genenum}, we can go
through the proof backwards and ask which Pfaffians correspond to the
plane partitions in Theorem~\ref{genenum} for general $c_1, c_2$. 
Of course, then we know that they evaluate to the nice factored
expression of Theorem~\ref{genenum}.

By exactly the method of \cite{min}, 
we can find a bijection to families of non-intersecting lattice paths
shown in Figure~\ref{paths}, express this number as a sum of
determinants using the main theorem on non-intersecting lattice paths
(see Lindstr\"om, \cite[Lemma~1]{LindAA} or Gessel and Viennot
\cite[Theorem~1]{gv}) ) and finally express this as a Pfaffian using a
theorem by Ishikawa and Wakayama (\cite{IshW95}).
\end{proof}

\begin{figure}
\begin{center}
\leavevmode
\epsfbox{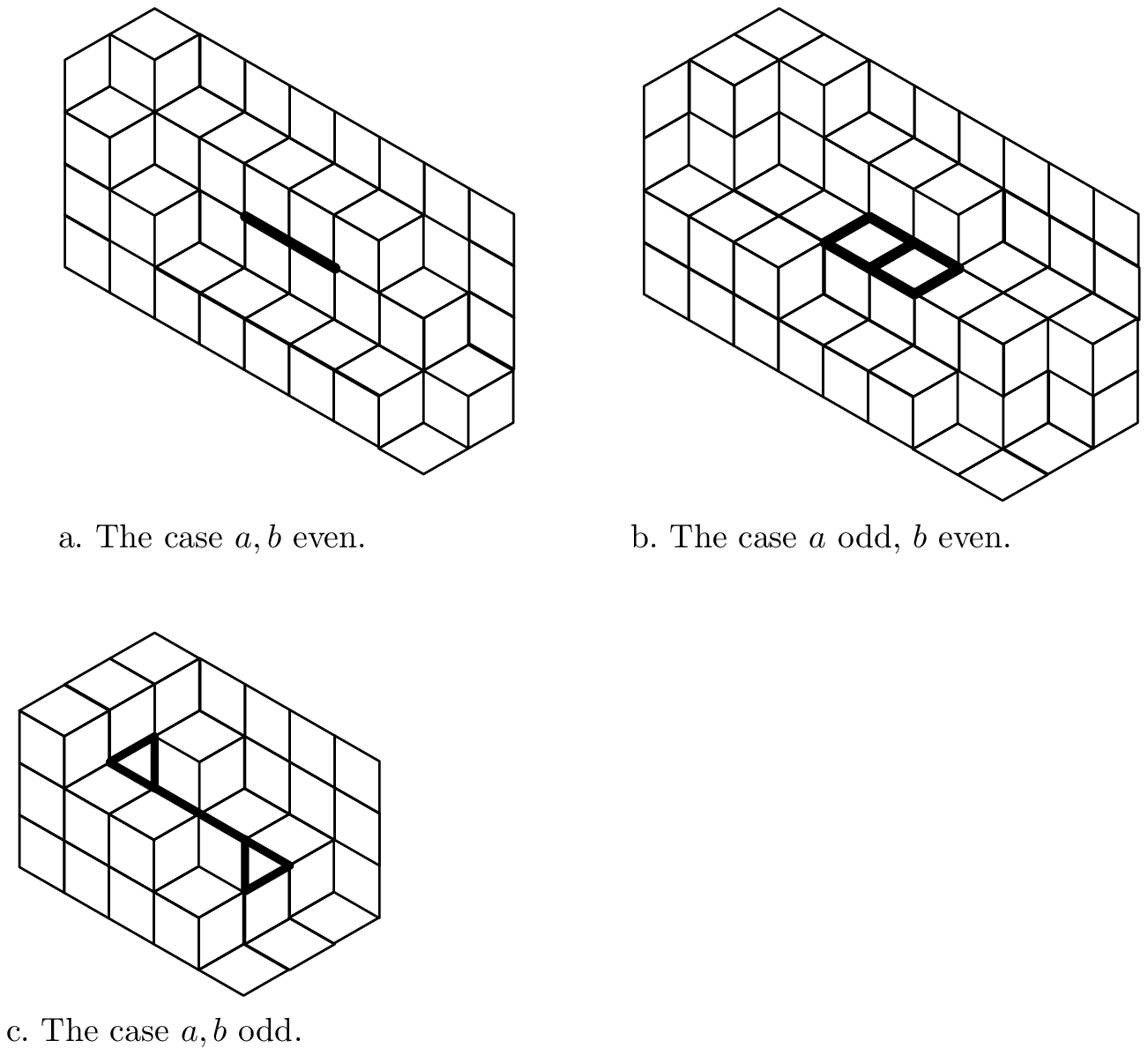}
\end{center}
\label{cuteee} 
\label{cuteoo} 
\label{cutoee} 
\label{cut}
\caption{Self-complementary plane partitions with a fixed line in the middle.}
\end{figure}

\begin{figure}
\begin{center}
\leavevmode
\epsfbox{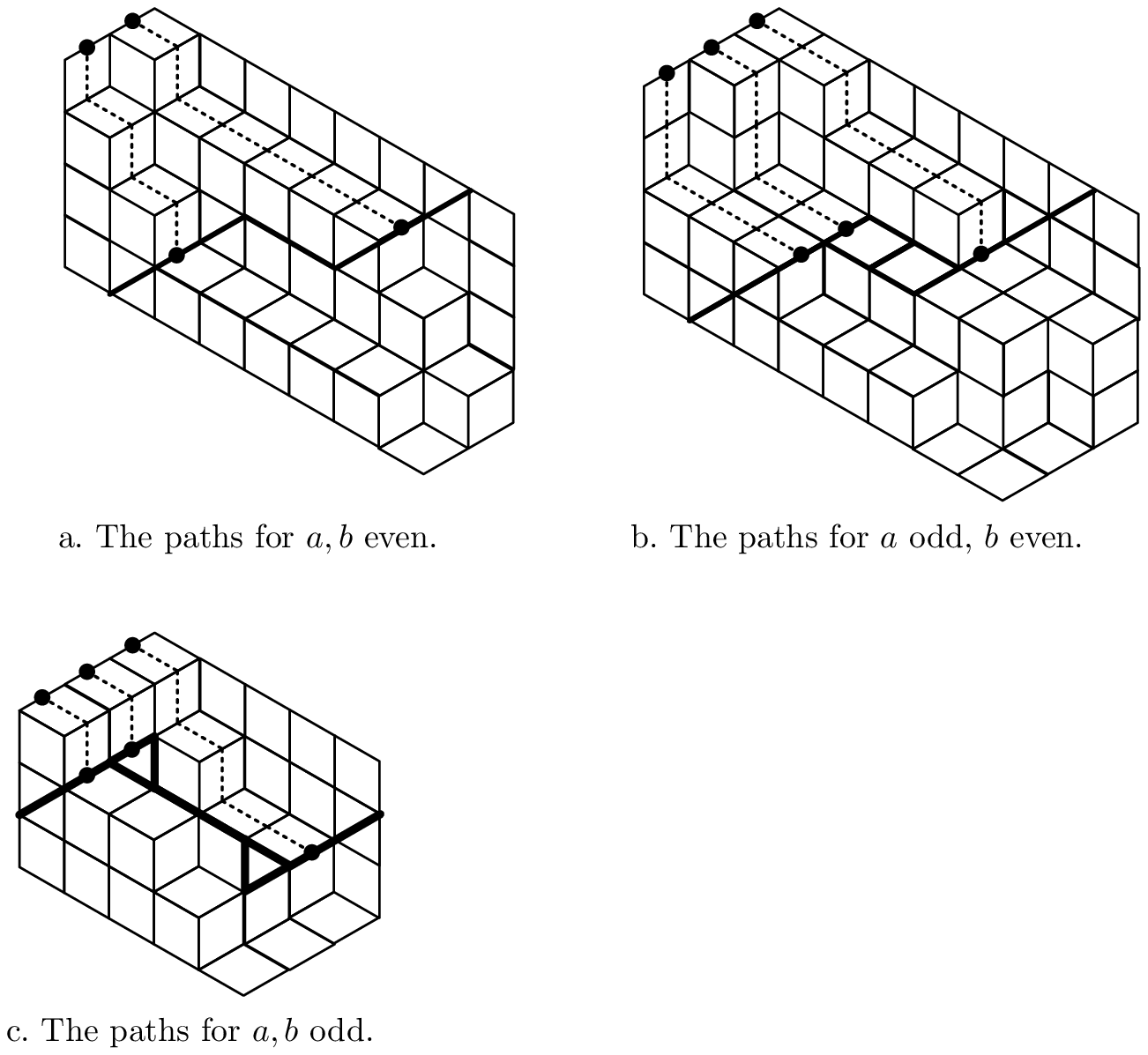}
\end{center}
\caption{The path families.}
\label{cuteeepaths} 
\label{cuteoopaths} 
\label{cutoeepaths} 
\label{paths}
\end{figure} 

\end{section}


\begin{thebibliography}{10}

\bibitem{An1}
G. E. Andrews, {\em Plane partitions (II): 
The equivalence of the Bender--Knuth and the MacMahon conjectures},
Pacific J. Math., {\bf 72} (1977), 283--291.

\bibitem{An2}
G. E. Andrews, {\em Plane partitions (I): 
The MacMahon conjecture},
Adv.\ in Math.\ Suppl.\ Studies, {\bf 1} (1978), 131--150.

\bibitem{An3}
G. E. Andrews, {\em Plane partitions V: The t.s.s.c.p.p.\ 
conjecture},  J.~Combin.\ Theory Ser.\ A, {\bf 66} (1994), 28--39.


\bibitem{min}
T.~Eisenk\"olbl, {\em $(-1)$-enumeration of plane partitions with
complementation symmetry},
Adv. in Appl. Math.  {\bf 30}, (2003),  no. 1-2, 53--95,
arXiv:math.CO/0011175. 

\bibitem{minsc}
T.~Eisenk\"olbl, {\em $(-1)$-enumeration of self-complementary plane
partitions}, arXiv:math.CO/0412118.


\bibitem{gv}
I.M.~Gessel, X.~Viennot, {\em Determinant, paths and plane
partitions}, Preprint, (1989).


\bibitem{IshW95}
M.~Ishikawa and M.~Wakayama, {\em Minor summation formula of Pfaffians},
Linear and Multilinear Algebra {\bf 39} (1995), 285--305.

\bibitem{IsOkTaZe}
M.~Ishikawa, S.~Okada, H.~Tagawa and J.~Zeng, {\em Generalizations of
  Cauchy's determinant and Schur's Pfaffian},
Adv. Appl. Math. {\bf 36} (2006), 251--287. {\tt arXiv:math.CO/0411280.}

\bibitem{Ku}
G.~Kuperberg, {\em An exploration of the permanent-determinant
method}, Electron.\  J.\  Combin. {\bf 5}, (1998), \#R46.
{\tt arXiv:math.CO/9810091.}

\bibitem{Ku2}
G. Kuperberg, {\em Symmetries of plane partitions and the 
permanent determinant method}, J.~Combin. Theory Ser. A {\bf 68}
(1994), 115--151.


\bibitem{LindAA}
B.    Lindstr\"om, {\em On the vector representations of induced
matroids}, Bull.\  London Math.\  Soc.\  {\bf 5} (1973), 85--90.

\bibitem{LR}
D.E. Littlewood and A. R. Richardson, {\em Group characters and
  algebra}, Philos. Trans. R. Soc., A, {\bf 233}, (1934), 99--141.

\bibitem{MM}
P.A.~MacMahon, {\em Combinatory Analysis\/}, vol.~2, Cambridge University
Press, (1916); reprinted by Chelsea, New York, (1960).

\bibitem{MRR2} W. H. Mills, D. P. Robbins and H. Rumsey,
{\em Proof of the Macdonald conjecture}, Inventiones Math. {\bf 66},
(1982), 73--87.

\bibitem{MRR} W. H. Mills, D. P. Robbins and H. Rumsey.
{\em Alternating sign matrices and descending plane partitions},
J. Combin.Theroy Ser. A {\bf 34} (1983), 340 -- 359.

\bibitem{Okad98} 
S.~Okada, {\em Applications of minor summation formulas to rectangular-shaped representations of classical
groups}, J.~Algebra, {\bf 205}, (1998), 337--367.

\bibitem{Stan3}
R.P. Stanley, {\em Theory and applications of plane partitions: Part
  2}, Stud. Appl. Math {\bf 50}, (1971), 259--279.

\bibitem{Stan1}
R. P. Stanley, {\em Enumerative Combinatorics}, Vol.~1,
Wadsworth \& Brooks/Cole, Pacific Grove, California, (1986).

\bibitem{Stan2} 
R.P. Stanley, {\em Symmetries of plane partitions},
J. Combin. Theory Ser\  A {\bf 43} (1986), 103--113;
Erratum {\bf 44} (1987), 310.



\bibitem{Stem94a}
J.R. Stembridge, {\em Some hidden relations involving the ten symmetry
classes of plane partitions},
J.~Combin.\   Theory Ser.\   A {\bf 68} (1994), 372--409. 

\bibitem{Stem94b}
J.R. Stembridge, {\em On minuscule representations, plane partitions
  and involutions in complex Lie groups},
Duke Math.\   J. {\bf 73} (1994), 469--490. 

\bibitem{Stem4}
J. R. Stembridge, {\em The enumeration of totally symmetric plane
partitions}, Adv. in Math. {\bf 111} (1995), 227--245.


\end{thebibliography}
\end{document}